\documentclass{article}
\usepackage{amsmath,amsthm,amsfonts,amssymb,url}
\usepackage{graphicx}
\title{Yuri Safarov (1958--2015)}
\date{}
\begin{document}
\maketitle
\begin{center}
\includegraphics[scale=0.3]{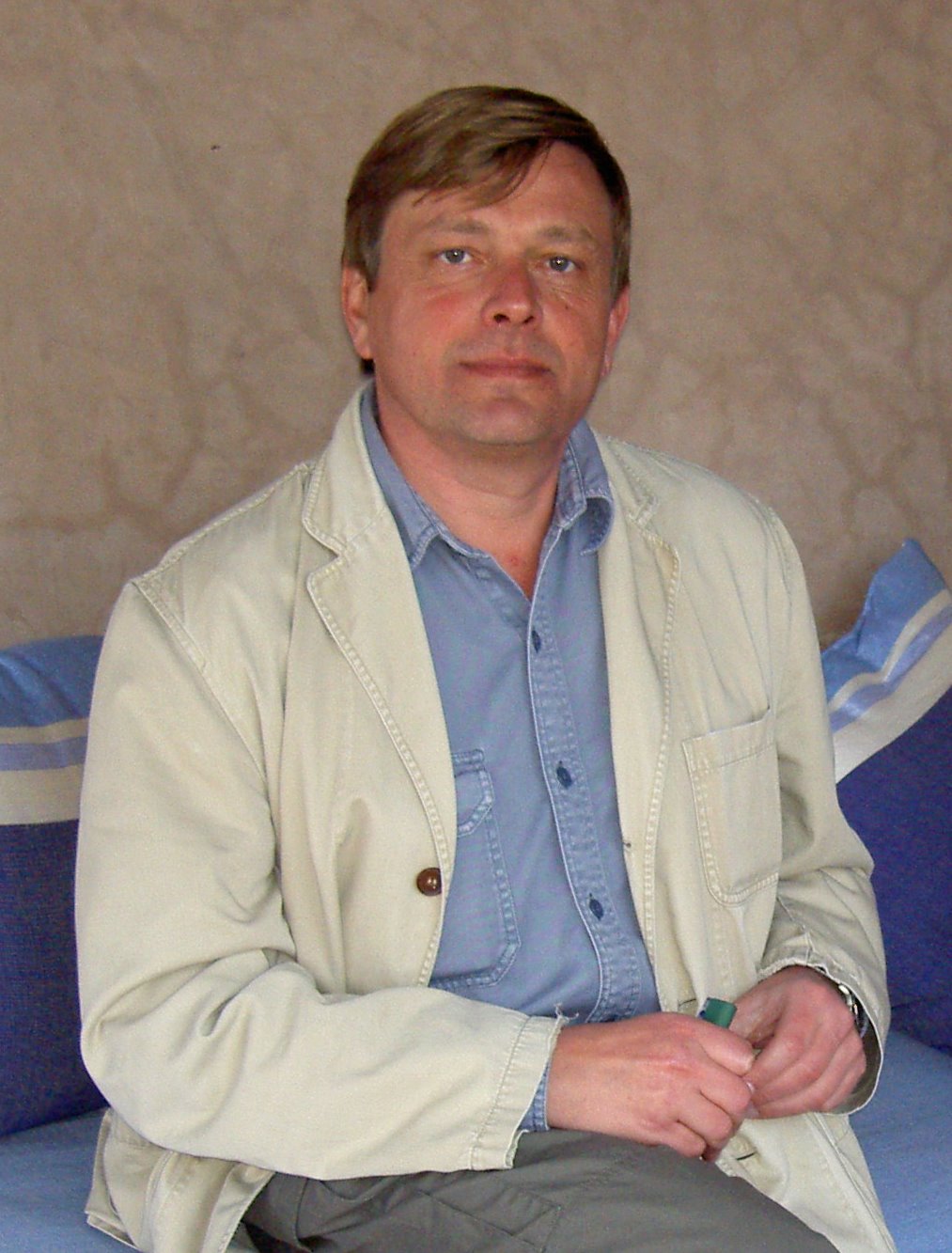}
\end{center}

\

\

This issue is dedicated to Yuri Safarov who passed away on 2
June 2015. It is a collection of papers in subject areas related to
Safarov's work which are written by people who knew him well.

The editor of this volume and author of this preface is a collaborator of Yuri Safarov.
The two of us coauthored the book \cite{mybook}
which contains Safarov's most influential mathematical results.
Lars H\"ormander, who reviewed the book for the Bulletin of the LMS,
summarised his review \cite{hormanders_review} as follows:
``In the reviewer's opinion, this book is indispensable
for serious students of spectral asymptotics''.

Yuri was born in Leningrad in the Soviet Union on 23 January 1958.
He spent his formative years in Leningrad (later renamed to Saint Petersburg),
before moving to London in 1993.
Yuri did his undergraduate studies at Leningrad State University, from
which he graduated in 1980. Upon graduation,
he was appointed Research Fellow at the
Leningrad Branch of the Steklov Mathematical Institute.

Yuri did a PhD at the Steklov Institute under the supervision of Michael Solomyak.
In his PhD thesis \cite{SafarovPhD} he considered two problems.

Problem 1: spectral problem for the operator
\begin{equation}
\label{maxwell operator}
\begin{pmatrix}
0&-i\operatorname{curl}\\
i\operatorname{curl}&0
\end{pmatrix}
\end{equation}
acting on pairs
$
\begin{pmatrix}
u_1\\
u_2
\end{pmatrix}
$
of solenoidal vector fields in a bounded domain $M\subset\mathbb{R}^3$
with smooth boundary $\partial M$.
The boundary conditions are
\begin{equation}
\label{maxwell boundary conditions}
\left.u_2\right|_{\partial M}=0,
\qquad
\left.u_1\times n\right|_{\partial M}=0,
\end{equation}
where $n$ is the normal to $\partial M$.
Safarov succeeded in overcoming two fundamental difficulties
evident in  \eqref{maxwell operator}, \eqref{maxwell boundary conditions}:
a) the difficulty that the problem is not elliptic and
b) the difficulty that the problem is not semi-bounded
(the eigenvalues $\lambda_k$ accumulate at $+\infty$ and $-\infty$).
Let us define the two counting functions
$N_\pm:(0,+\infty)\to\mathbb{N}\cup\{0\}$
as
$\displaystyle N_\pm(\lambda):=
\sum_{0<\pm\lambda_k<\lambda}1\,.$
The function $N_+(\lambda)$ counts eigenvalues $\lambda_k$
between zero and $\lambda$,
whereas
the function $N_-(\lambda)$ counts eigenvalues $\lambda_k$
between $-\lambda$ and zero.
Safarov proved that
\[
N_\pm(\lambda)=\frac{V}{3\pi^2}\lambda^3+O(\lambda^2)
\qquad\text{as}\qquad
\lambda\to+\infty,
\]
where $V$ is the volume of the domain $M$.
Furthermore, under the nonperiodicity condition
(not too many periodic billiard trajectories, see \cite{mybook} for details)
Safarov proved the stronger result
\[
N_\pm(\lambda)=\frac{V}{3\pi^2}\lambda^3+o(\lambda^2)
\qquad\text{as}\qquad
\lambda\to+\infty.
\]

Problem 2: transmission problem.
Here we have two compact $n$-dimensional
Riemannian manifolds, $M_+$ and $M_-$,
connected by a common smooth boundary $\Gamma$ and the aim is to
study the spectrum of the Laplacian on the manifold $M:=M_+\cup M_-$
subject to natural transmission boundary conditions on the $(n-1)$-dimensional
hypersurface $\Gamma$. Here the fundamental difficulty,
one which Safarov succeeded in overcoming, is that the billiard flow experiences
branching on the hypersurface $\Gamma$:
in addition to the usual reflected trajectory one gets
a refracted trajectory.
Under appropriate assumptions on the branching billiard flow
(nonperiodicity plus the condition that there are not too many trajectories
experiencing an infinite number of reflections in a finite time), Safarov obtained a
two-term asymptotic formula for the counting function of this problem.
In fact, Safarov went further and examined the situation when there are many periodic
trajectories in which case he derived an asymptotic formula with a periodic
factor in the second term.

In 1990 the Steklov Institute awarded Yuri a DSc (habilitation).
In the Soviet Union it was most unusual to be awarded a DSc at such a young age.
Yuri's DSc thesis \cite{SafarovDSc}
contained a number of important results, some of which
are listed below.

Firstly, Safarov performed a comprehensive analysis of the situation
when there are many periodic trajectories. For a particular type of
spectral problem such an analysis already started in \cite{SafarovPhD}.
However, in \cite{SafarovDSc} it was done in far greater generality. Also,
Safarov showed how to deal with the (very common) case when the
periodic factor in the second asymptotic term of the counting function is discontinuous.
Non-classical (i.e.~not purely polynomial) two-term asymptotic formulae
for the counting function are Yuri's trademark results.

Secondly, Safarov established very precise asymptotic formulae for the spectral function.
The spectral function is a ``local'' version of the counting function: whereas the
counting function depends only on one variable, the spectral parameter $\lambda$,
the spectral function depends also on the point $y\in M$ or even a pair of
points $x,y\in M$ (here $M$ is the manifold on which the operator acts).
For the case of a scalar positive definite operator the spectral function is defined below,
see formula \eqref{definition of spectral function}.
Without going into technicalities let us just mention that establishing asymptotic
formulae for the spectral function is a delicate matter and the geometric
conditions involved in the derivation of such formulae are different from those
for the ``global'' counting function.
Asymptotic formulae for the spectral function are also
Yuri's trademark results.

Thirdly, Safarov examined elliptic  systems on manifolds without boundary.
Let $A$ be an elliptic (i.e.~determinant of principal symbol is nonzero)
self-adjoint first order linear pseudodifferential operator acting on $m$-columns
of complex-valued half-densities over a connected
closed (i.e.~compact and without boundary) $n$-dimensional  manifold $M$.
We also assume that the eigenvalues of the principal
symbol are simple, but we do not assume that the operator is semi-bounded.
Under the above assumptions Safarov derived an explicit formula
for the principal symbols of the $m$ individual oscillatory integrals
making up the propagator $e^{-itA}$, $t\in\mathbb{R}$ being the time variable.
Safarov's formula was subsequently heavily used by a number of authors, myself included.
The formula is nontrivial in that it differs fundamentally from the well known
scalar formula of Duistermaat and Guillemin.

I met Yuri around 1983 and,
despite the fact that we lived in different cities (Yuri in Leningrad and I in Moscow),
we started interacting and exchanging ideas.
Towards the end of the 1980s Simon Gindikin suggested that Yuri and I write
a book for the American Mathematical Society. We took up the offer but our work
progressed very slowly because we were preoccupied with difficulties of everyday
life in precollapse USSR and because we did not have access to personal computers
which would have allowed us to typeset the text in TeX.
But things changed when the two of us moved to the UK:
I took up an appointment at Sussex in 1991 and Yuri one at King's College London
in 1993. Yuri was initially appointed as an SERC Advanced Fellow.
He became a Professor in 1997 and continued
working at King's until his untimely death.

We finished writing our book \cite{mybook} in 1996.
The book deals with an elliptic self-adjoint scalar differential operator of even order $2m$
acting on half-densities over
a compact $n$-dimensional manifold $M$ with smooth boundary $\partial M$,
subject to $m$ boundary conditions on $\partial M$.
Without loss of generality, the operator is assumed to be positive definite
and the spectral parameter is denoted by $\lambda^{2m}$, which simplifies
subsequent asymptotic formulae.
The objects of study are the counting function
\begin{equation}
\label{definition of global counting function}
N(\lambda):=
\sum_{\lambda_k<\lambda}1
\end{equation}
and the spectral function
\begin{equation}
\label{definition of spectral function}
e(\lambda,x,y):=
\sum_{\lambda_k<\lambda}v_k(x)\,\overline{v_k(y)}\,,
\end{equation}
where the $\lambda_k$ and the $v_k(x)$, $k=1,2,\ldots$,
are the eigenvalues and the orthonormalised
eigenfunctions. Obviously, the counting function
\eqref{definition of global counting function}
is the trace
of the spectral projection, whereas the spectral function
\eqref{definition of spectral function}
is the integral kernel of the spectral projection.

For the counting function $N(\lambda)$ the book establishes explicit two-term asymptotic
formulae as $\lambda\to+\infty$. Here the second term may be classical
(constant times $\lambda^{n-1}$) or non-classical
(periodic function of $\lambda$ times~$\lambda^{n-1}$).

The spectral function $e(\lambda,x,y)$ is examined in the book
on the diagonal and in the interior of the manifold,
i.e.~under the assumptions $x=y\not\in\partial M$.
It is shown that the classical second term in the asymptotics of
$e(\lambda,y,y)$ is zero, whereas the non-classical one may be nonzero.
Safarov provided a simple example of a situation when
the second term in the asymptotics of
$e(\lambda,y,y)$ is non-classical: this is the example of the Laplacian
in a planar domain bounded by an ellipse, subject to the Dirichlet or Neumann boundary
condition. When the point $y$ is a focus of the ellipse
we get a non-classical second term which Safarov wrote out explicitly.

The book uses a version of Levitan's hyperbolic equation method,
so a propagator had to be constructed using microlocal techniques.
Yuri insisted on using his particular microlocal technique
based on an invariantly defined complex-valued phase function.
This technique originates from the paper
\cite{complex_phase_paper}
and allows one to construct oscillatory integrals globally in time,
passing without problems through caustics and avoiding the need for tackling
compositions of oscillatory integrals.
Furthermore, this approach leads to a new, purely analytic and straightforward,
definition of the Maslov index.

In the latter part of his career, after finishing the book
 \cite{mybook}, Yuri studied a variety of different questions
in the spectral theory of differential operators.
Some of these publications are mentioned below.

In \cite{szego} Laptev and Safarov established Szeg\H o type limit theorems
for elliptic semi-bounded self-adjoint pseudodifferential operators.

In \cite{SafarovPLMS1997} Safarov developed
a coordinate-free calculus of pseudodifferential opera\-tors
on a manifold. Here the critical idea, originally due to Widom, is the use of a connection.

In \cite{SafarovJFA2001} Safarov performed a comprehensive study
of Fourier Tauberian theorems,
coming up with an `ultimate',  very general version of a Fourier Tauberian theorem.
He used it for the study of
the counting function of the Dirichlet Laplacian in an arbitrary  domain 
of finite volume.

In \cite{NetrusovSafarov} Netrusov and Safarov
studied the asymptotic distribution of eigen\-values of the Laplacian on a
bounded domain in $\mathbb{R}^n$ with rough boundary.
The authors obtained an explicit remainder estimate in the
Weyl formula for the Dirichlet Laplacian.
Furthermore, they established
sufficient conditions for the validity of the Weyl formula for the Neumann Laplacian.

In \cite{Lapointe} Lapointe, Polterovich and Safarov
studied the average growth of the spectral function on a Riemannian manifold.
Two types of averaging were considered: with respect to the
spectral parameter and with respect to a point on a manifold.

In \cite{Filonov} Filonov and Safarov studied, in the abstract setting, the following question:
if the self-commutator $[A^*,A]$ of an operator $A$ is small, is $A$ close to a normal (or
even a diagonal) operator?

In \cite{McKeag} Safarov, jointly with his PhD student McKeag,
revisited the issue of coordinate-free description
of pseudodifferential operators on manifolds. This paper was published at a time when
I also became interested in pseudodifferential operators and Fourier integral operators
on manifolds, and we had a series of fruitful discussions with Yuri in the last few years of
his life. These discussions had a profound influence on my current research.

In \cite{Jakobson} Jakobson, Safarov and Strohmaier studied
the semiclassical limit of spectral theory on manifolds whose metrics
have jump-like discontinuities
and proved a quantum ergodicity theorem for such systems.

Safarov's last paper \cite{Kachkovskiy}, written jointly with his PhD student Kachkovskiy,
was published posthumously. It revisits the self-commutator problem
which was initially exam\-ined in \cite{Filonov}.
The problem at hand is one posed by Halmos \cite{halmos}
and can be informally stated as follows: how close is an almost
normal matrix to a normal one?
Kachkovskiy's and Safarov's paper
\cite{Kachkovskiy}
effectively settles this long-standing problem.

Safarov's contributions to mathematics were recognised by mathematical
societies both in the USSR and the UK.
The Leningrad Mathematical Society awarded Yuri the 1987 Young Mathematician Prize
and 
The London Mathematical Society awarded him the 1996  Junior Whitehead Prize.
The citation \cite{LMS} for the latter read as follows.

\emph{
Dr Y. Safarov of King's College, London, is awarded a Junior Whitehead Prize
for his work on the spectral analysis of partial differential operators and microlocal
analysis.}

\emph{
He undertook in the late eighties a detailed investigation of compact manifolds for
which the spectral asymptotics of the Laplace--Beltrami operator has a non-classical
second term. The complete analysis of this involved looking at the geometrical
structure of the set of periodic orbits and their relationship with the behaviour of the
wave equation --- a physically-motivated approach which required completely new
insights to make into a piece of rigorous mathematics.}

\emph{
Since this result, Dr Safarov has been engaged in several other major pieces of
work, including a detailed asymptotic analysis of the spectral function, and a
powerful abstract version of Szeg\H o's theorem. He has also suggested a new approach
to the theory of Fourier integral operators based on the use of a complex, rather than
real, phase function. Amongst the advantages of this viewpoint is a new definition of
the Maslov index.}

\emph{
Most recently, Safarov has developed a coordinate-free approach to the theory of
pseudo-differential operators, introducing Weyl symbols for a pseudo-differential
operator on a manifold with connection.}

Throughout his career
Yuri supervised many PhD students and several postdocs. He also volunteered
to teach a course \emph{Distributions, Fourier Transforms and Microlocal Analysis}
at the London Taught Course Centre, an organisation which provides taught lecture
courses for PhD students from several London universities and some universities from
surrounding areas.
He continued lecturing
at the London Taught Course Centre
when he was already ill.
Yuri was an excellent lecturer
at all university teaching levels, undergraduate and postgraduate.

Despite the progressing illness,
Yuri actively supervised his last PhD student Paolo Battistotti
and signed off the final version of the thesis just days before his death.

Outside of mathematics Yuri's passion was
chess. As a youngster, participating in a tournament
he played against the future world
champion Garry Kasparov and that game
ended in a draw. Even though Kasparov was
five years younger than Yuri, this was a remarkable
achievement.
In 1981 Yuri was awarded the title Master of Sport of the USSR.
He also played for the chess team of Leningrad.
At some point in his life Yuri had to decide
whether to pursue a professional chess career
or become a professional mathematician. Fortunately
for the subject of spectral theory he
chose the latter.

Yuri was an open, friendly and outgoing person who enjoyed life to the fullest.
He will be greatly missed by friends and colleagues.

\

\

\

\hfill Dmitri Vassiliev

\end{document}